\documentclass[12pt]{amsart}
\usepackage{amsmath,amssymb}
\usepackage{color}
\theoremstyle{plain}
\newtheorem{thm}{Theorem}[section]
\newtheorem{lem}[thm]{Lemma}

\newtheorem{pro}[thm]{Proposition}

\begin{document}
\title{Homogeneous spaces not separated by arcs}

\author{Alexandre Karassev}
\address{Department of Computer Science and Mathematics, Nipissing University,
100 College Drive, P.O. Box 5002, North Bay, ON, P1B 8L7, Canada}
\email{alexandk@nipissingu.ca}
\thanks{The first author was partially supported by NSERC Grant 257231-15}

\author{Vesko  Valov}
\address{Department of Computer Science and Mathematics, Nipissing University,
100 College Drive, P.O. Box 5002, North Bay, ON, P1B 8L7, Canada}
\email{veskov@nipissingu.ca}
\thanks{The second author was partially supported by NSERC Grant 261914-19}

\begin{abstract}
It was shown in
\cite{vmv} that regions in strongly locally homogeneous locally compact metric spaces of dimension $\ge 2$ are not separated by arcs. We improve this result by replacing strong local homogeneity with homogeneity. Moreover, we prove the result for the case when only one end point of an arc is in the interior of the region.
\end{abstract}

\makeatletter
\@namedef{subjclassname@2020}{\textup{2020} Mathematics Subject Classification}
\makeatother



\maketitle

\section{Introduction}

All spaces are assumed to be metric separable, and all maps continuous. Everywhere below, $\mathbb S^1$ denotes the circle. Recall that a space $X$ is called:

\begin{itemize}

\item[-] {\it homogeneous} if for any $x,y\in X$ there exists a homeomorphism $f\colon X\to Y$ such that $f(x)=y.$

\medskip

\item[-] {\it strongly locally homogeneous} (SLH) if for any $x\in X$ there exists a base of neighbourhoods $\mathcal B$ of $x$ such that for any $B\in \mathcal B$ and any $y,z\in B$ there exists a homeomorphism $f\colon X\to X$ such that $f(y)=z$ and $f$ is the identity on $X\setminus B$.

\end{itemize}

Connected SLH spaces are homogeneous but the converse is not the case (e.g. any solenoid is homogeneous but not SLH). It was shown in \cite{KP} that locally connected, simply connected, homogeneous metric continuum cannot be separated by arcs. The question of whether the assumption of simply connectedness is necessary was also posed in \cite{KP}. In \cite{vmv} it was shown that for any locally compact SLH space $X$ of dimension $\ge 2$ any its region $G$ cannot be separated by an arc lying entirely in $G$. Here by a region in a space $X$ we mean an open connected subspace of $X$. In this paper we improve the result of \cite{vmv} and provide a positive answer to the question from \cite{KP} as follows.

\begin{thm}\label{main} Let $X$ be a locally compact homogeneous space with ${\rm dim }\, X = 2$ and $G$ be a region in $X$. Let $J$ be an arc in $X$ such that  $J\setminus\{p\}\subset G$ where $p$ is an endpoint of $J$. Then $J$ does not separate $G$.
\end{thm}

Note that results of \cite{kru} imply that no region of homogeneous locally compact space $X$ with ${\rm dim} \, X = n\ge 1$ can be separated by a closed subset of dimension $\le n-2$. Thus Theorem~\ref{main} also holds in case ${\rm dim}\, X>2.$

\section{Preliminary statements}

Suppose $f\colon A\to \mathbb S^1$ is a map defined on a closed subset $A$ of a compact space $K$ that cannot be extended over $K$. Applying Zorn's lemma, one can show that there  exists a compactum $M=M(A,f)\subset K$ containing $A$ such that $f$ cannot be extended over $M$, but it is extendable over any closed proper subset of $M$ containing $A$ (note that $M(A,f)$ is not unique). Such a set $M(A,f)$ is called an {\em $f$-membrane spanned on $A$}.

We also use the following notation. Suppose $A$ is partition in a space $Z$ between two closed disjoint sets $P, Q\subset Z$. Then there exist two open disjoint subset $W_P$ and  $W_Q$ of $Z$ containing $P$ and $Q$, respectively, such that $Z\backslash A=W_P\cup W_Q$. We denote
$\Lambda_P=W_P\cup A$ and $\Lambda_Q=W_Q\cup A$.

\begin{lem}\label{membrane} Let $X$ be as in Theorem~\ref{main}.
 For every $x\in X$ there exists a compactum $\Phi$ and a local base $\mathcal B_x$ in $X$ with the following property.
For every $U,V\in\mathcal B_x$ with $\overline V\subset U$ there exist a map
$f_{U,V}:\Phi(U,V)\to\mathbb S^1$, where $\Phi(U,V)=(\overline U\backslash V)\cap \Phi$ such that:
\begin{itemize}
\item[(i)] There exists a compactum $M_V\subset\overline{V}$ which is a membrane spanned on $\Phi\cap{\rm Bd}\, \overline V$ for the map $f_V=f_{U,V}|(\Phi\cap {\rm Bd}\, \overline V)$;
\item[(ii)] $f_U=f_{U,V}|(\Phi\cap{\rm Bd}\, \overline  U)$ is not extendable over $\Phi\cap\overline U$;
\item[(iii)] If $A$ is a partition in $\Phi(U,V)$ between $\Phi\cap{\rm Bd}\, \overline U$ and $\Phi\cap{\rm Bd}\overline V$, then
$f_A=f_{U,V}|A$ is extendable over $\Lambda_{{\rm Bd}\, \overline V}\cup B$ for any proper closed set $B\subset M_V$, but not extendable over $\Lambda_{{\rm Bd}\, \overline V}\cup\overline V$.
\end{itemize}
\end{lem}

\begin{proof}
Since $X$ is a countable union of compact sets, there exists a compactum $Y\subset X$ with $\dim Y=2$ (otherwise, by the countable sum theorem for $\dim$, $\dim X\leq 1$). Since $\dim Y=2$ there exists a proper closed subset $F\subset Y$ and  a map $f\colon F\to\mathbb S^1$ such that $f$ is not extendable over $Y$. Consequently, there exists a
compact set $\Phi\subset Y$ containing $F$ such that $\Phi$ is an $f$-membrane spanned on $F$. Since $X$ is homogeneous, we may assume that $x\in \Phi\backslash F$.
Now, let $\mathcal B_x$ be the family of all
open sets $U$ in $X$ such that:
$$U={\rm int}\, \overline U \mbox{ and } \overline U\cap F=\varnothing.$$
Suppose $U,V\in\mathcal B_x$ with $\overline V\subset U$. Then
$\Phi\backslash V$ is a proper closed subset of $\Phi$ containing $F$. Hence there exists $f'\colon \Phi\backslash V\to\mathbb S^1$ extending $f$ such that $f'$ is not extendable over $\Phi\cap\overline V$.
Let $f_{U,V}=f'|\Phi(U,V)$.
Then $f_V=f_{U,V}|(\Phi\cap{\rm Bd}\, \overline V)$ is not extendable over $\Phi\cap\overline V$, otherwise $f$ would be extendable over $\Phi$. Therefore there exists a compactum $M_V\subset\overline V$ which is an $f_V$-membrane spanned on $\Phi\cap {\rm Bd}\, \overline V$.
Similarly,
$f_U=f_{U,V}|(\Phi\cap {\rm Bd}\, \overline U)$ is not extendable over $\Phi\cap\overline U$.

To prove item $(iii)$, suppose $A$ is a partition in $\Phi(U,V)$ between $\Phi\cap {\rm Bd}\,\overline{U}$ and $\Phi\cap{\rm Bd}\, \overline{V}$. Then $f_A$ is the restriction of $f_{\Lambda_{{\rm Bd}\, \overline V}}$ on $A$, where
$f_{\Lambda_{{\rm Bd}\, \overline V}}=f_{U,V}|\Lambda_{{\rm Bd}\, \overline V}$. Moreover, $f_V$ is extendable over any proper closed set $B\subset M_V$.
Hence, $f_A$ is extendable over $\Lambda_{{\rm Bd}\, \overline V}\cup B$ for any closed proper set $B\subset M_V$. On the other hand, $f_A$ is the restriction of
$f_{\Lambda_{{\rm Bd}\, \overline U}}=f_{U,V}|\Lambda_{{\rm Bd}\, \overline U}$ on $A$. This implies that $f_A$
is not extendable over $\Lambda_{{\rm Bd}\, \overline V}\cup\overline V$ because so is $f_{\Lambda_{{\rm Bd}\, \overline U}}$.
\end{proof}

We need the following version of Effros' theorem \cite{e} for locally compact spaces, see \cite{kru}.

\begin{thm}\label{effros}
Let $X$ be a homogeneous locally compact space with a metric $\rho$, $a\in X$ and $\varepsilon>0$. Then there exists $\delta>0$ such that for every $x\in X$ with $\rho(x,a)<\delta$ there exists a homeomorphism $h\colon X\to X$ with $h(a)=x$ and $\rho(h(y),y)<\varepsilon$ for all $y\in X$.
\end{thm}

The following is a slight modification of Proposition~2.3 from \cite{vmv}, the proof of which is identical to that given in \cite{vmv}.

\begin{pro}\label{nullhomotopic}
Let	$\varphi\colon A\to\mathbb S^1$ be a map defined on a closed subspace of a space $Y$. Suppose that closed subsets $Y_1$ and $Y_2$ of $Y$ satisfy the following conditions:

\begin{itemize}
	\item $Y=Y_1\cup Y_2;$
	\item $\left|Y_1\cap Y_2\cap A\right| \le 1;$
	\item $\varphi|Y_i\cap A$ is extendable over $Y_i$, $i=1,2;$
	\item $\varphi$ is not extendable over $Y.$
	
\end{itemize}

Then there exists $\psi\colon Y_1\cap Y_2\to \mathbb S^1$ which is not nullhomotopic.
\end{pro}

\section{Main theorem}

In this section we prove Theorem~\ref{main}.
Let $X$ be a homogeneous locally compact separable metric space and $G$ be a region in $X$ with ${\rm dim}\, X = 2$. Note that the homogeneity of $X$ and the countable sum theorem implies that $X$ is everywhere $2$-dimensional. Let $J$ be an arc in $G$ such that one of its endpoints, say $b$, is in $G$, and the other endpoint, say $a$, is in ${\rm Bd}\, G$, and $J\setminus \{a\}\subset G$ (the proof for the case when $J\subset G$ is analogous). We will show that $J$ does not separate $G$. Assume it does. As in the proof of Lemma~3.2 from \cite{vmv}, we may assume that there are two disjoint non-empty open subsets of $G$ with $G\setminus J = G_1\cup G_2$, and $\overline{G_1}\cap\overline{G_2}\cap G = S\setminus \{a\}$, such that  $b=\max\{s\colon s\in S\}$ and $S$ is a closed subset of $J$ without isolated points.

By Lemma~\ref{membrane} there exist a compactum $\Phi$ containing the point $b$ and a local base $\mathcal B_b$ of $b$ in $X$ satisfying the statement of that lemma.
We may assume that $\overline U\subset G$ for all $U\in\mathcal B_b$. Following  notations of  Lemma~\ref{membrane}, we may also assume that $b\in M_V\setminus{\rm Bd}\, V$, where $U,V\in\mathcal B_b$ with $\overline V\subset U$. Indeed,
applying Theorem~\ref{effros}, we can find $\delta>0$ corresponding to the point $b$ and  $\varepsilon={\rm dist}\,(\overline U,X\backslash G)$, and choose $V$ so small that its diameter is less than $\delta.$
Then there exists a homeomorphism $h\colon X\to X$ which is $\delta$-close to identity and such that $h(b)\in M_V\setminus{\rm Bd}\, V$. Now we can consider the sets $h^{-1}(U),h^{-1}(V),$ and $h^{-1}(M_V)$ instead of $U,V,$ and $M_V$, respectively.

Further, since $f_V$ cannot be extended over $M_V$, we have $\dim M_V=2$. This implies that either $M_V\cap G_1\neq\varnothing$ or  $M_V\cap G_2\neq\varnothing$. If $M_V\cap G_1\neq\varnothing$, we take $x\in M_V\cap G_1$ and let $\varepsilon_1=\min\{\alpha,\beta\}$, where
$\alpha={\rm dist}({\rm Bd}\, U,X\setminus G)$ and $\beta={\rm dist}(x,\overline G_2\cup{\rm Bd}\, V)$. Applying again Theorem~\ref{effros}, we can push $M_V$ towards $G_2$ by a $\delta_1$-small homeomorphism $h$,  where $\delta_1$ corresponds to $b$ and $\varepsilon_1,$ such that $h(M_V)\cap G_i\neq\varnothing$, $i=1,2$, and $h(V)$ still contains $b$. Thus, everywhere below we can assume that $b\in M_V\setminus {\rm Bd}\, V$ and
$M_V\setminus {\rm Bd}\,V$ meets both $G_1$ and $G_2.$

Since $S\cap{\rm Bd}\,V$ and $S\cap{\rm Bd}\,U$ are closed disjoint subsets of $S\cap \Phi(U,V)$ there exists a partition $H=\{c\}\in S$  in
$S\cap \Phi(U,V)$ between $S\cap{\rm Bd}\,V$ and $S\cap{\rm Bd}\,U$ (if $S\cap(U\setminus\overline V)\cap \Phi=\varnothing$, we take $H=\varnothing$).
By \cite[Corollary 3.5]{vm} there exists a partition $T$ in $\Phi(U,V)$ between $\Phi\cap{\rm Bd}\,V$ and $\Phi\cap{\rm Bd}\,U$ with $T\cap S\subset H$.
Then $A=T\cup H$ is a partition in $\Phi(U,V)$ between $\Phi\cap{\rm Bd}\,V$ and $\Phi\cap{\rm Bd}\,U$ with $A\cap S=H$.
Let $\widetilde M=M_V\cup\Lambda_{{\rm Bd}\,\overline V}$. There are two possible cases:
\begin{itemize}
\item[(1)] $f_A$ is not extendable over $\widetilde M$;
\item[(2)] $f_A$ is extendable over $\widetilde M$.
\end{itemize}
We will show that each of these cases is  impossible. Indeed, assume first that  $f_A$ is not extendable over $\widetilde M,$ and let
$S_U=S\cap \overline U$ and $K=\widetilde M\cup S_U$. Since $f_A$ is not extendable over $\widetilde M$, $f_A$ is not extendable over
$K$. Let $P_i=S_U\cup (K\cap\overline G_i)$, $i=1,2$. Obviously, $P_1\cup P_2=K$ and $P_1\cap P_2=S_U$.
Since each $M_V\cap\overline G_i$ is a proper subset of $M_V$, $f_A$ is extendable over each of the sets $A\cup(K\cap\overline G_i)$.
Finally, since $\dim S_U\leq 1$, $f_A$ is extendable over each of the sets $A\cup P_i$. Therefore, we can apply Proposition~\ref{nullhomotopic} and conclude that there exists a map $\beta\colon S\to\mathbb S^1$ which is not null-homotopic. This is a contradiction because $S\subset J$ and every map from $J$ to $\mathbb S^1$ is null-homotopic.

Now suppose that $f_A$ is extendable over $\widetilde M$.
Since $f_A$ is not extendable over $\Lambda_{{\rm bd}\overline V}\cup\overline V$, there exists a minimal subset $M_A$ of $\Lambda_{{\rm bd}\overline V}\cup\overline V$ containing $\widetilde M$ such that $M_A$ is an $f_A$-membrane spanned on $A$. Since $M_A$ contains $M_V$, $M_A$ meets both $G_1$ and $G_2$.
Hence, we can proceed as in the previous case with $M_V$ replaced by $M_A,$ and obtain again a contradiction.

\end{document}